\documentstyle[12pt]{amsart}

\newcommand{\nz}{\mbox{${\Bbb Z}$}}

\newcommand{\vep}{\varepsilon}
\newcommand{\vphi}{\varphi}
\theoremstyle{plain}
 \newtheorem{thm}{Theorem}[section]
 \newtheorem{prop}[thm]{Proposition}

\theoremstyle{definition}
 \newtheorem{defn}{Definition}[section]
 
 \newtheorem{ex}{Example}[section]
\theoremstyle{remark}

 \newtheorem{ack}{Acknowledgment}
      
\begin{document}
\title{Generalized Drinfeld realization of  quantum superalgebras and 
$U_q(\hat {\frak osp}(1,2))$}
\author{Jintai Ding}
\author{Boris  Feigin}
\address{Jintai Ding, Department of Mathematical Sciences, University of Cincinnati}
\address{Boris Feigin, Landau Institute of Theoretical Physics}
\maketitle
\centerline{Dedicated to our friend Moshe Flato} 
\begin{abstract}
In this paper, we extend the generalization of 
Drinfeld realization of quantum affine algebras to  quantum affine superalgebras
with its Drinfeld comultiplication and its Hopf algebra structure, 
which depends on a function $g(z)$ satisfying the relation: 
$$g(z)=g(z^{-1})^{-1}.$$
In particular, we present the Drinfeld realization of $U_q(\hat {\frak osp}(1,2))$
and its Serre relations. 

\end{abstract}

\pagestyle{plain}
\section{Introduction.} 

Quantum groups  as a noncommutative 
and noncocommutative  Hopf algebras were   discovered by Drinfeld\cite{Dr1} 
and Jimbo\cite{J1}. 
The standard  definition of a  quantum group is given  as a deformation 
of universal enveloping algebra of a simple (super-)Lie algebra
 by the basic generators and the relations based on the data coming from the  
corresponding Cartan matrix. However, for the case of quantum affine algebras, 
there is a different aspect of the theory, namely their loop 
realizations. 
The first approach was given by Faddeev, Reshetikhin and
Takhtajan \cite{FRT} and Reshetikhin and Semenov-Tian-Shansky
\cite{RS},  who obtained a realization  of
the quantum loop algebra $U_q(
{\frak{ g}} \otimes {\Bbb C}[t, t^{-1}])$ via a canonical 
 solution
 of the Yang-Baxter equation depending on a parameter $z\in\Bbb C$. 
 On the other hand, Drinfeld \cite{Dr2} gave another realization  
of the quantum affine algebra $U_q(\hat {\frak g})$ and its special degeneration 
called the Yangian, 
which  is widely  used in constructions of special representation of
affine quantum algebras\cite{FJ}.  
In \cite{Dr2}, Drinfeld only gave the realization of 
the quantum affine algebras as an algebra, and as an  algebra this 
realization is equivalent to  the approach above \cite{DF} through
certain Gauss decomposition for the case of $U_q(\hat {\frak gl}(n))$. 
Certainly, the most important aspect of the structures 
 of the quantum groups is its Hopf algebra
structure, especially its comultiplication. Drinfeld also 
constructed a new Hopf algebra structure for this loop realization.
The new  comultiplication in this formulation, which we call  the  Drinfeld 
comultiplication, is simple and has very important applications \cite{DM} \cite{DI2}.

In \cite{DI}, we observe that in the Drinfeld realization  of 
quantum affine algebras $U_q(\hat {\frak sl}(n))$, the structure constants are certain 
rational functions $g_{ij}(z)$, whose 
 functional property of $g_{ij}(z)$ decides completely the Hopf algebra structure. 
In particular, for the case of $U_q(\hat {\frak sl}_2)$, its Drinfeld realization
is given completely in terms of a function $g(z)$, which has the following function property:
$$g(z)=g(z^{-1})^{-1}.$$
This leads  us to 
generalize this type of Hopf algebras. Namely, we can substitute  
$g_{ij}(z)$ by other functions that satisfy the functional property 
of $g_{ij}(z)$, to derive new Hopf algebras. 

In this paper, we will further  extend the generalization of
the Drinfeld realization of  $U_q(\hat {\frak sl}_2)$ to derive quantum affine superalgebras. 
As an example, we will also present the quantum affine superalgebra 
$U_q(\hat {\frak osp}(1,2))$ in terms of the new formulation, in particular, we present 
the Serre relations in terms of the current operators. 

The paper is organized as the following: in Section 2, we recall the main results in 
\cite{DI} about the generalization of Drinfeld realization of $U_q(\hat {\frak sl}(2))$; 
in Section 3, we present the definition of the generalized Drinfeld realization of 
quantum superalgebras; in Section 4, we present the formulation of 
$U_q(\hat {\frak osp}(1,2)).$

\section{}
In \cite{DI}, we derive a generalization of Drinfeld realization of 
$U_q(\hat {\frak sl}_n)$. For the case of 
$U_q({\frak sl}_ 2)$, we first present the complete definition. 

Let $g(z)$ be an 
  analytic functions  satisfying the 
following property that $g(z)=g(z^{-1})^{-1}$ and
$\delta(z)$ be  the distribution with support at $1$. 

\begin{defn}
$U_q(g,f{\frak sl}_2)$ is an
 associative algebra with unit 
1 and the generators:  $x^{\pm}(z)$, 
$\vphi (z)$,
$\psi (z)$, a central element $c$
 and a nonzero  complex parameter $q$, where $z\in {\Bbb C}^*$. 
$\vphi (z)$ and $\psi (z)$ are invertible. 
In terms of the generating functions:
the defining relations are 
\begin{align*}
& \vphi (z)\vphi (w)=\vphi (w)\vphi(z), \\
& \psi(z)\psi  (w)=\psi  (w)\psi(z), \\
& \vphi(z)\psi  (w)\vphi(z)^{-1}\psi  (w)^{-1}=
  \frac{g(z/wq^{-c})}{g(z/wq^{c})}, \\
& \vphi(z)x  ^{\pm}(w)\vphi(z)^{-1}=
  g(z/wq^{\mp \frac{1}{2}c})^{\pm1}x  ^{\pm}(w), \\
& \psi(z)x  ^{\pm}(w)\psi(z)^{-1}=
  g(w/zq^{\mp \frac{1}{2}c})^{\mp1}x  ^{\pm}(w), \\
& [x^+(z),x  ^-(w)]=\frac{1}{q-q^{-1}}
  \left( \delta(\frac zwq^{-c})\psi(wq^{\frac{1}{2}c})-
          \delta(\frac zwq^{c})\vphi(zq^{\frac{1}{2}c}) \right), \\
&  x^{\pm}(z)x  ^{\pm}(w)=
  g (z/w)^{\pm 1}x ^{\pm}(w)x^{\pm}(z). 
\end{align*}
\end{defn}

\begin{thm}
The algebra $U_q(g,f{\frak sl}_2)$  
has a Hopf algebra structure, which are given 
by the following formulae. 

\noindent{\bf Coproduct $\Delta$}
\begin{align*}
\text{(0)}& \quad \Delta(q^c)=q^c\otimes q^c, \\
\text{(1)}& \quad \Delta(x   ^+(z))=x   ^+(z)\otimes 1+
            \vphi   (zq^{\frac{c_1}{2}})\otimes x   ^+(zq^{c_1}), \\
\text{(2)}& \quad \Delta(x   ^-(z))=1\otimes x   ^-(z)+
            x   ^-(zq^{c_2})\otimes \psi   (zq^{\frac{c_2}{2}}), \\
\text{(3)}& \quad \Delta(\vphi   (z))=
            \vphi   (zq^{-\frac{c_2}{2}})\otimes\vphi   (zq^{\frac{c_1}{2}}), \\
\text{(4)}& \quad \Delta(\psi   (z))=
            \psi   (zq^{\frac{c_2}{2}})\otimes\psi   (zq^{-\frac{c_1}{2}}),
\end{align*}
where $c_1=c\otimes 1$ and $c_2=1\otimes c$.

\noindent{\bf Counit $\vep$}
\begin{align*}
\vep(q^c)=1 & \quad \vep(\vphi   (z))=\vep(\psi   (z))=1, \\
            & \quad \vep(x   ^{\pm}(z))=0.
\end{align*}
\noindent{\bf Antipode $\quad a$}
\begin{align*}
\text{(0)}& \quad a(q^c)=q^{-c}, \\
\text{(1)}& \quad a(x   ^+(z))=-\vphi   (zq^{-\frac{c}{2}})^{-1}
                               x   ^+(zq^{-c}), \\
\text{(2)}& \quad a(x   ^-(z))=-x   ^-(zq^{-c})
                               \psi   (zq^{-\frac{c}{2}})^{-1}, \\
\text{(3)}& \quad a(\vphi   (z))=\vphi   (z)^{-1}, \\
\text{(4)}& \quad a(\psi   (z))=\psi   (z)^{-1}.
\end{align*}

\end{thm}

Strictly speaking, $U_q(g,f{\frak sl}_2)$ is not an algebra. 
This concept, which  we call a functional algebra, has  already 
been used before \cite{S}, etc.

The Drinfeld realization for the case of 
$U_q(\hat {\frak sl}_2)$ \cite{Dr2} as a Hopf algebra is  different, and 
it an algebra and Hopf algebra
defined with current operators in terms of formal power series.

Let $g(z)$ be an analytic functions that satisfying the 
following property that $g(z)=g(z^{-1})^{-1}=
G^+(z)/G^-(z)$, where $G^\pm (z)$ is an analytic function
 without poles except at $0$ or $\infty$ and $G^\pm(z)$ have no common
zero point.
Let $ \delta(z)=\sum_{n\in \Bbb Z}z^n$, where $z$ is a formal variable.

\begin{defn}
The algebra $U_q(g,{\frak sl}_2)$ is an associative algebra with unit 
1 and the generators: $\bar a   (l)$,$\bar b   (l)$, $x^{\pm}   (l)$, for 
 $l\in \nz $ and a central element $c$. 
Let $z$ be a formal variable and 
$$x   ^{\pm}(z)=\sum_{l\in \nz}x   ^{\pm}(l)z^{-l},$$
$$\vphi   (z)=\sum_{m\in \nz}\vphi   (m)z^{-m}=\exp [
\sum_{m\in \nz_{\leq 0}}\bar a   (m)z^{-m}]\exp [
\sum_{m\in \nz_{> 0}}\bar a   (m)z^{-m}]$$
  and 
$$\psi   (z)=\sum_{m\in \nz}\psi   (m)z^{-m}=\exp [
\sum_{m\in \nz_{\leq 0}}\bar b   (m)z^{-m}]\exp [
\sum_{m\in \nz_{>0}}\bar b   (m)z^{-m}]. $$
 In terms of the 
formal variables $z$, $w$, 
the defining relations are 
\begin{align*}
& a   (l)a  (m)=a  (m)a   (l), \\
& b   (l)b  (m)=b  (m)b   (l), \\
& \vphi   (z)\psi  (w)\vphi   (z)^{-1}\psi  (w)^{-1}=
  \frac{g (z/wq^{-c})}{g (z/wq^{c})}, \\
& \vphi   (z)x  ^{\pm}(w)\vphi   (z)^{-1}=
  g (z/wq^{\mp \frac{1}{2}c})^{\pm1}x  ^{\pm}(w), \\
& \psi   (z)x  ^{\pm}(w)\psi   (z)^{-1}=
  g (w/zq^{\mp \frac{1}{2}c})^{\mp1}x  ^{\pm}(w), \\
& [x   ^+(z),x  ^-(w)]=\frac{1}{q-q^{-1}}
  \left( \delta(\frac zwq^{-c})\psi   (wq^{\frac{1}{2}c})-
          \delta(\frac zwq^{c})\vphi   (zq^{\frac{1}{2}c}) \right), \\
& G^{\mp}(z/w) x   ^{\pm}(z)x  ^{\pm}(w)=
  G^{\pm }(z/w)x  ^{\pm}(w)x   ^{\pm}(z), 
\end{align*}
where
  by  $g(z)$ we   mean the Laurent expansion of $g(z)$ 
in a region $r_{1}>|z|> r_{2}$. 
\end{defn}

\begin{thm}
The algebra $U_q(g,{\frak sl}_2)$  
has a Hopf algebra structure.  The formulas for the coproduct $\Delta$, 
the counit $\vep$ and the antipode $a$ are the same as given in 
Theorem 2.1. 
\end{thm}

Here, one has to be careful with the 
expansion of the structure functions $g(z)$ and 
$\delta(z)$, for the  reason that the  relations between $x^{\pm}   (z)$
and $x^{\pm}  (z)$ are different from the case of the  functional algebra
above. 

\begin{ex}

Let $\bar g(z)$ be a an analytic function such that $\bar g(z^{-1})=-z^{-1}\bar g(z)$. Let 
$g(z)=q^{-2}\frac{\bar g(q^{2}z)}{\bar g(q^{-2}z)}$. Then 
$g(z)=g(z^{-1})^{-1}$. 
With this  $g(z)$, we can define
an algebra $U_q(g,{\frak sl}_2)$ by  Definition 2.2.

{\bf Case I.} 
Let $\bar g(z)$=$1-z$. 
$\vphi   (m)=\psi   (-m)=0$ for $m\in \nz_{>0}$ and $\vphi   (0)\psi   (0)=1$. 
Then this 
algebra is  $U_q(\hat {\frak sl}_2)$. 
\end{ex}

{\bf Case II}
Let $\bar g(z)= \theta_p(z)=\prod_{j>0}(1-p^j)(1-p^{j-1}z)(1-p^jz^{-1})$ be the Jacobi's
theta function.
$\theta_p(z^{-1})=-z^{-1}\theta_p(z)=\theta_p(pz)$. 
 We will call this 
algebra $U_q(\theta, {\frak sl}_2)$. We take the expansion of $g(z)$
in the region $|q^2|>|z|>|q^2p|$.
If we take the limit that 
 $p$ goes to zero, we would 
recover $U_q(\hat {\frak sl}_2)$.

\section{}
Similarly, we can use the same idea to define quantum affine superalgebras, which 
is  an  extension of the generalization  in the section above. 
Again, we will start from the same  function $g(z)$, namely,  
let $g(z)$ be an 
  analytic functions  satisfying the 
following property that $g(z)=g(z^{-1})^{-1}$ and
$\delta(z)$ be  the distribution with support at $1$.

\begin{defn}
$U_q(g,f{\frak s})$ is an
 a ${\Bbb  Z}_2$ graded associative algebra with unit 
1 and the generators:  $x^{\pm}(z)$, 
$\vphi (z)$,
$\psi (z)$, a central element $c$
 and a nonzero  complex parameter $q$, where $z\in {\Bbb C}^*$, 
$x^{\pm}(z)$ are  graded 1(mod2), and $\vphi (z)$,
$\psi (z)$ and c are graded 0(mod2). 
$\vphi (z)$ and $\psi (z)$ are invertible. 
In terms of the generating functions:
the defining relations are 
\begin{align*}
& \vphi (z)\vphi (w)=\vphi (w)\vphi(z), \\
& \psi(z)\psi  (w)=\psi  (w)\psi(z), \\
& \vphi(z)\psi  (w)\vphi(z)^{-1}\psi  (w)^{-1}=
  \frac{g(z/wq^{-c})}{g(z/wq^{c})}, \\
& \vphi(z)x  ^{\pm}(w)\vphi(z)^{-1}=
  g(z/wq^{\mp \frac{1}{2}c})^{\pm1}x  ^{\pm}(w), \\
& \psi(z)x  ^{\pm}(w)\psi(z)^{-1}=
  g(w/zq^{\mp \frac{1}{2}c})^{\mp1}x  ^{\pm}(w), \\
& \{x^+(z),x  ^-(w)\}=\frac{1}{q-q^{-1}}
  \left( \delta(\frac zwq^{-c})\psi(wq^{\frac{1}{2}c})-
          \delta(\frac zwq^{c})\vphi(zq^{\frac{1}{2}c}) \right), \\
&  x^{\pm}(z)x  ^{\pm}(w)=-
  g (z/w)^{\pm 1}x ^{\pm}(w)x^{\pm}(z), 
\end{align*}
where $\{ x, y\}=xy+yx$. 
\end{defn}

 Here we remark that the above relations are basically the same as in that of Definition 2.1
except the relation between $x^{\pm}(z)$ and $x  ^{\pm}(w)$ respectively, which differs by 
a negative sign. 

Accordingly  we have that,  for the tensor algebra, the multiplication is defined for 
 homogeneous elements $a,$  $b,$ $c $, $d$ by
$$
(a\otimes b)(c\otimes d)=(-1)^{[b][c]}\,(ac\otimes bd),
$$
where $[a] \in{\bf Z}_2$ 
denotes the grading of the element $a$.

Similarly we have:

\begin{thm}
The algebra $U_q(g,f{\frak s})$  
has a  graded Hopf algebra structure, whose coproduct, counit
and antipode  are given 
by the same formulae of $U_q(q, f{\frak sl}_2)$ in Theorem 2.1. 
\end{thm}

As for the case of $U_q(g,f{\frak sl}_2)$ is not a graded  algebra but rather a 
a graded functional algebra. 

Let 
 $$g(z)=g(z^{-1})^{-1}=
G^+(z)/G^-(z),$$ where $G^\pm (z)$ is an analytic function
 without poles except at $0$ or $\infty$ and $G^\pm(z)$ have no common
zero point.

\begin{defn}
The algebra $U_q(g,{\frak s})$ is 
 ${\Bbb  Z}_2$ graded associative algebra with unit 
1 and the generators: $\bar a   (l)$,$\bar b   (l)$, $x^{\pm}   (l)$, for 
 $l\in \nz $ and a central element $c$, where 
$x^{\pm}(l)$ are graded 1(mod 2) and the rest are graded o(mod 2). 
Let $z$ be a formal variable and 
$$x   ^{\pm}(z)=\sum_{l\in \nz}x   ^{\pm}(l)z^{-l},$$
$$\vphi   (z)=\sum_{m\in \nz}\vphi   (m)z^{-m}=\exp [
\sum_{m\in \nz_{\leq 0}}\bar a   (m)z^{-m}]\exp [
\sum_{m\in \nz_{> 0}}\bar a   (m)z^{-m}]$$
  and 
$$\psi   (z)=\sum_{m\in \nz}\psi   (m)z^{-m}=\exp [
\sum_{m\in \nz_{\leq 0}}\bar b   (m)z^{-m}]\exp [
\sum_{m\in \nz_{>0}}\bar b   (m)z^{-m}]. $$
 In terms of the 
formal variables $z$, $w$, 
the defining relations are 
\begin{align*}
& \vphi (z)\vphi (w)=\vphi (w)\vphi(z), \\
& \psi(z)\psi  (w)=\psi  (w)\psi(z), \\
& \vphi(z)\psi  (w)\vphi(z)^{-1}\psi  (w)^{-1}=
  \frac{g(z/wq^{-c})}{g(z/wq^{c})}, \\
& \vphi(z)x  ^{\pm}(w)\vphi(z)^{-1}=
  g(z/wq^{\mp \frac{1}{2}c})^{\pm1}x  ^{\pm}(w), \\
& \psi(z)x  ^{\pm}(w)\psi(z)^{-1}=
  g(w/zq^{\mp \frac{1}{2}c})^{\mp1}x  ^{\pm}(w), \\
& \{x^+(z),x  ^-(w)\}=\frac{1}{q-q^{-1}}
  \left( \delta(\frac z wq^{-c})\psi(wq^{\frac{1}{2}c})-
          \delta(\frac zwq^{c})\vphi(zq^{\frac{1}{2}c}) \right), \\
&  (G^\mp (z/w))^{}x^{\pm}(z)x  ^{\pm}(w)=-
  (G^\pm (z/w))^{}x ^{\pm}(w)x^{\pm}(z),
\end{align*}
where
  by  $g(z)$ we   mean the Laurent expansion of $g(z)$ 
in a region $r_{1}>|z|> r_{2}$. 
\end{defn}

The above relations are basically the same as in that of Definition 2.2
except the relation between $x^{\pm}(z)$ and $x  ^{\pm}(w)$ respectively, which  differs by 
a negative sign. 
The  
expansion direction of the structure functions $g(z)$ and 
$\delta(z)$ is very important, for the  reason that the  relations between $x^{\pm}   (z)$
and $x^{\pm}  (z)$ are different from the case of the  functional algebra
above.

\begin{thm}
The algebra $U_q(g,{\frak s})$  
has a Hopf algebra structure.  The formulas for the coproduct $\Delta$, 
the counit $\vep$ and the antipode $a$ are the same as given in 
Theorem 2.1. 
\end{thm}

\begin{ex}

Let $\bar g(z)= 1$.
From  \cite{CJWW}\cite{Z},  we can see that 
$U_q(1, {\frak s})$ is basically the same as $U_q(\hat {\frak gl}(1,1))$.

\end{ex}

\section{}

For a rational  function $g(z)$ that satisfies 
$$ g(z)=g(z^{-1})^{-1},$$
it is clear that $g(z)$ is determined by its poles and its zeros, which are paired to 
satisfy the relations above.
For the simplest case  (except $g(z)=1$) 
 that $g(z)$ has only one pole and one zero, we have 
$$ g(z)= \frac {zp-1}{z-p}, $$
where $p$ is the location of the pole of $g(z)$. 
Unfortunately, in this case, we do not know anyway to identify the algebra 
$U_q(g, {\frak s})$ with other know structures. 

Then, comes the second simplest case that  
$g(z)$ has two poles and two zeros.  
We know also that $g(z)$ must be in the form  that 
$$ g(z)= \frac {zp_1-1}{z-p_1}\frac {zp_2-1}{z-p_2}. $$ 
In this section, we shall establish that for this case, 
$U_q(g, {\frak s})$  is related to the affine quantum superalgebra
$U_q(\hat{\frak osp}(1,2))$. 

\vskip .5in 
{}From now on, let us fix $g(z)$ to be 
$\frac {zp_1-1}{z-p_1}\frac {zp_2-1}{z-p_2}$. 
\vskip .5in 

As in \cite{DM}\cite{DK}, for the case of quantum affine algebras, it is 
very important to understand the poles and zero of the product of current 
operators. We will start with the relations between $X^+(z)$ with itself. 
{}From the definition, we know that 

$$(z-p_1w)(z-p_2w) X^+(z)X^+(w)= - (zp_1-w)(zp_2-w) X^+(w)X^+(z). $$

{}From this, we know that 
$ X^+(z)X^+(w)$ has two poles, which are located at $(z-p_1w)=0$ and $(z-p_2w)=0$. 

This also implies that 
\begin{prop}
 $ X^+(z)X^+(w)=0$, when $z=w$. 
\end{prop}

If we assume that $U_q(g,{\frak s})$ is related to some quantized affine superalgebra, 
then we can see that the best chance we have is $U_q(\hat {\frak osp}(1,2))$ by looking 
at the number of zeros and poles of $ X^+(z)X^+(w)$.

However, for the case of $U_q(\hat {\frak osp}(1,2))$, we know we need an extra Serre relation. 
For this,  we will follow  the idea in \cite{FO}. 

Let 
Let $$f(z_1, z_2)= {(z_1-p_1z_2)(z_1-p_2z_2)}.$$ 
$$Y^+(z,w)= \frac {(z-p_1w)(z-p_2w)} {z-w} X^+(z)X^+(w), $$
$$Y^+(z_1,z_2,z_3)= \frac {(z_1-p_1z_2)(z_1-p_2z_2)} {z_1-z_2}
\frac {(z_1-p_1z_3)(z_1-p_2z_3)} {z_1-z_3}\times $$
$$ \frac {(z_2-p_1z_3)(z_2-p_2z_3)} {z_2-z_3}
 X^+(z_1)X^+(z_2)X^+(z_3),$$
$$F(z_1,z_2,z_3)=  {(z_1-p_1z_2)(z_1-p_2z_2)}
{(z_3-p_1z_1)(z_3-p_2z_1)} 
 {(z_2-p_1z_3)(z_2-p_2z_3)}= $$
$$f(z_1,z_2)f(z_2,z_3)f(z_3,z_1),$$
$$\bar F(z_1,z_2,z_3)= 
f(z_2,z_1)f(z_2,z_3)f(z_3,z_1). $$

Let $V_(z_1,z_2,z_3)$ be the  algebraic variety of the zeros of $F(z_1,z_2,z_3)$. 
Let  $ V(a(z_1),a(z_2),a(z_3)),$  be the image of the action of $a$ on this variety, 
where  $a\in$ $S_3$,  the permutation group on $z_1,z_2,z_3$. 
Let $\bar V_(z_1,z_2,z_3)$ be the  algebraic variety of the zeros of $\bar F(z_1,z_2,z_3)$. 
Let  $\bar  V(a(z_1),a(z_2),a(z_3)),$  be the image of the action of $a$ on this variety, 
where  $a\in$ $S_3$.

\begin{prop} 
$Y^+(z,w)$ has no poles and is symmetry with respect to $z$ and $w$. 
$Y^+(z_1,z_2,z_3)$ has no poles and is symmetry with respect to $z_1$, $z_2$, $z_3$. 
\end{prop} 

Following  the idea in \cite{FO}, we would like to define  the following  conditions
that may be imposed on  our algebra. 

{\bf Zero Condition I:}

$Y^+(z_1,z_2,z_3)$ is zero on  at least one line  that crosses 
(0,0,0),  and   this line must lie in a 
 $ V(a(z_1),a(z_2),a(z_3))$ for some element $a\in S_3$ 
\vskip .5in

{\bf Zero Condition II:}
$Y^+(z_1,z_2,z_3)$ is zero on  at least one line  that crosses 
(0,0,0),  and   this line must lie in a 
 $ \bar V(a(z_1),a(z_2),a(z_3))$ for some element $a\in S_3$ 
\vskip .5in 

For the line that crosses $(0,0,0)$ where  $Y^+(z_1,z_2,z_3)$ is zero, we call it the zero line of 
$Y^+(z_1,z_2,z_3)$. 
Because $Y^+(z_1,z_2,z_3)$ is symmetric with respect to the action of $S_3$ on 
$z_1,z_2,z_3$, if a line is the zero line of $Y^+(z_1,z_2,z_3)$, then clearly 
the orbit of the line under the action of $S_3$ is also an  zero line.

\vskip .15in

{\bf  Remark 1.} 
There is a simple symmetry that we would prefer  to  choose the function $F(z_1, z_2, z_3)$ 
to determine the 
variety $V(z_1, z_2 , z_3)$. 
We have that $$F(z_1, z_2, z_3)= f(z_1,z_2)f(z_2,z_3)f(z_3,z_1).$$  
Let $S_2^1$ be the permutation group acting on $z_1,z_2.$
Let $S_2^2$ be the permutation group acting on $z_2,z_3.$
Let $S_2^1$ be the permutation group acting on $z_3,z_1.$
Clearly, \cite{FO}
 we can choose from a family of varieties determined by the functions 
$ f(a_1(z_1),a_1(z_2))f(a_2(z_2),a_2(z_3))f(a_3(z_3),a_3(z_1))$ for 
$a_1\in S_2^1$, $a_2\in S_2^2$, $a_3\in S_2^3$.   
For each  such a function $$f(a_1(z_1),a_1(z_2))f(a_2(z_2),a_2(z_3))f(a_3(z_3),a_3(z_1)),$$ 
we can attach a oriented diagram, whose nods are $z_1,z_2, z_3$, and the arrows are 
given by  $ (a_1(z_1)\rightarrow a_1(z_2)),$  $(a_2(z_2)\rightarrow a_2(z_3))$ $ (a_3(z_3)\rightarrow a_3(z_1))$. For example the diagram of $F(z_1,z_2,z_3)$ is 
given by 

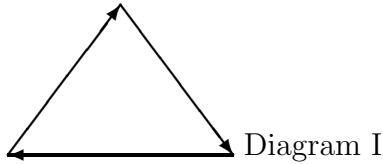
\begin{center}
\setlength{\unitlength}{.5 cm}
\begin{picture}(6,4)\thicklines  
\put(6 ,0  ){\vector (-1,0){6 }} 
\put(0  ,0  ){\vector (3,4){3 }}
 \put(3,4){\vector (3,-4){3}} 
\end{picture}
Diagram I
\end{center} 

While $f(z_2,z_1)f(z_2,z_3)f(z_3,z_1)$ is denoted by 

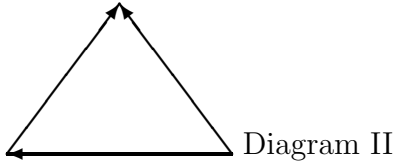
\begin{center}
\setlength{\unitlength}{.5 cm}
\begin{picture}(6,4)\thicklines  
\put(6 ,0  ){\vector (-1,0){6 }} 
\put(0  ,0  ){\vector (3,4){3 }}
 \put(6,0){\vector (-3,4){3}} 
\end{picture}
Diagram II
\end{center} 

It is not difficult to see that the diagram for 
$F(z_1,z_2,z_3)$ is symmetric in the sense that all the points are equivalent, but  for the 
second situation,  the top point $z_1$ is different from the other two, in the sense that 
there are two arrows coming to $z_1$, one to $z_3$ and none to $z_2$. 
There is only one other such a diagram given by 
$F(z_1, z_3, z_2)$, which however comes from the $S_3$ action on $F(z_1,z_2,z_3)$. 
\vskip .5in

{}From  the above, we have the following: 
\begin{prop}
Under the action of $S_3$ on the  the family of varieties  determined by the functions 
$ f(a_1(z_1),a_1(z_2))f(a_2(z_2),a_2(z_3))f(a_3(z_3),a_3(z_1))$ for 
$a_1\in S_2^1$, $a_2\in S_2^2$, $a_3\in S_2^3$, there are two orbits. One of the orbit consists
of  the two varieties determined by 
$F(z_1,z_2,z_3)$ and $F(z_1,z_2,z_3)$; and the rest forms another orbit.   
\end{prop} 

This shows that indeed we have two choices with respect the  zero conditions:
the Zero condition I and the Zero condition II.

\vskip .5in 

We know that a zero line 
is always in the form $$z_1=q_1z_2=q_2z_3. $$ 

Then we have 
\begin{prop} If we impose  the Zero condition I on the  algebra
$U_q(g,{\frak s})$, we have 
$$p_1=p_2^{-2}, $$ or
$$p_2=p_1^{-2}. $$   
\end{prop} 

{\bf Proof}. 
The proof is very simple. 
Because of the action of $S_3$, we know that one of the line 
must lie in $V(z_1,z_2,z_3)$. 
Let assume this line to be $z_1=z_2q_1=q_2z_3$. 

We know immediately that 
$q_1$ must be $p_1$ or $p_2$. 

Let us first deal with the case that 
$q_1=p_1$.  
We also know that $q_2$ must be either $p_1^{-1}$ or $p_2^{-1}$ by looking that the relations between 
$z_1$ and $z_3$. 

{\bf Case 1} Let $q_2=p_1^{-1}$, which implies that $p_1^{-2}$ must be either $p_1$ or $p_2$. 
Clearly, it can not be $p_1$, which implies that the impossible condition  $p_1=1$

Therefore, we have that 
$$p_1^{-2}=p_2,$$
which is what we want.

{\bf Case 2} Let $q_2=p_2^{-1}$, which implies that 
$p_1p_2$ must be either $p_1$ or $p_2$. 
Clearly, it can not be $p_1$ because it implies $p_2=1$, it can 
not be neither be $p_2$, which implies  that $p_1=1$. 

This completes the proof for   $$p_2=p_1^{-2}. $$ 

Similarly, if we  have that $q_1=p_2$, we can, then,  show
 $$p_1=p_2^{-2}. $$

However from the algebraic point of view, the two condition are 
equivalent in the sense that $p_1$ and $p_2$ are symmetric.

Also we have that 
\begin{prop}  If we impose the  Zero condition II on the  algebra $U_q(g,{\frak s})$, we have 
$$p_1=p_2^{2}, $$ or
$$p_2=p_1^{2}. $$   
\end{prop} 

However we also have that: 

{\bf If we impose the Zero condition II on the algebra $U_q(g,{\frak s})$, then 
 $U_q(g,{\frak s})$ is not  a Hopf algebra anymore. }

The reason is that the Zero condition II can not be satisfied 
by comultiplication, which can be checked by direct calculation.

This is the  most important reason that we will choose  
the Zero condition I to be imposed on the algebra 
$U_q(g,{\frak s})$, which  comes actually from the consideration of 
Hopf algebra structure. 
Namely, if   we choose  $V(z_1,z_2,z_3)$ or the  equivalent 
ones  which has the same diagram  presentation as Diagram I 
to define the zero line of  $Y^+(z_1,z_2,z_3)$, then, the quotient 
algebra derived from the 
{\bf Zero Condition I} is still a Hopf algebra with the same Hopf algebra structure (comultiplication, 
counit and antipode).

\vskip .5in
{}From now on, we impose the Zero condition I on the algebra 
$U_q(g,{\frak s})$, and   
 let us fix the notation  such that 
$$p_1=q^2, $$
$$p_2=q^{-1}. $$ 
\vskip .5in 

Similarly, we define

$$Y^-(z,w)= \frac {(z-p^{-1}_1w)(z-p^{-1}_2w)} {z-w} X^+(z)X^+(w), $$
$$Y^-(z_1,z_2,z_3)= \frac {(z_1-p^{-1}_1z_2)(z_1-p^{-1}_2z_2)} {z_1-z_2}
\frac {(z_1-p^{-1}_1z_3)(z_1-p^{-1}_2z_3)} {z_1-z_3} \times $$
$$\frac {(z_2-p^{-1}_1z_3)(z_2-p^{-1}_2z_3)} {z_2-z_3}
 X^+(z_1)X^+(z_2)X^+(z_3),$$

We now define the q-Serre relation. 

\vskip .5in 
{\bf q-Serre relations}
 
$Y^+(z_1,z_2,z_3)$ is zero on  the  line 
$$z_1=z_2q^{-1}=z_3q^{-2}. $$ 
$Y^-(z_1,z_2,z_3)$ is zero on  the  line 
$$z_1=z_2q^{}=z_3q^{2}. $$
\vskip .15in

The q-Serre relations can also be formulated in more algebraic way. 

\begin{prop}
The q-Serre relations are equivalent to the following two relations:
$$ 
 \frac {(z_3-z_1q^{-1})(z_3-z_1q^{3})(z_1-z_2q^2)}{z_3-z_1q}X^+(z_3)X^+(z_1)X^+(z_2) +$$
$$
\frac { (z_2-z_1q^2) ({z_2-z_1q^{-1}}) (z_3-z_1q^{-1})(z_3-z_1q^{3})}{
(z_1-z_2q^{-1}) (z_3-z_1q) }X^+(z_3) X^+(z_2)X^+(z_1) )  -$$
$$
  \frac{((z_1-z_3q^2)(z_1-z_3q)(z_1q-z_3q^{-1})(z_1-z_2q^2)} {(z_1-z_3)(z_3-z_1q^2)}
X^+(z_1)X^+(z_2)X^+(z_3)- $$
$$\frac{((z_1-z_3q^2)(z_1-z_3q)(z_1q-z_3q^{-1})(z_2-z_1q^2)(z_2-z_1q^{-1})}
 {(z_1-z_3)(z_3-z_1q^2)(z_1-z_2q^{-1})}
 X^+(z_2)X^+(z_1)X^+(z_3) )= 0 ,$$

$$ 
 \frac {(z_3-z_1q^{})(z_3-z_1q^{-3})(z_1-z_2q^{-2})}{z_3-z_1q^{-1}}X^-(z_3)X^-(z_1)X^-(z_2) + $$
$$\frac { (z_2-z_1q^{-2}) ({z_2-z_1q^{}}) (z_3-z_1q^{})(z_3-z_1q^{-3})}{
(z_1-z_2q^{}) (z_3-z_1q^{-1}) }X^-(z_3) X^-(z_2)X^-(z_1) )  -$$
$$
  \frac{((z_1-z_3q^{-2})(z_1-z_3q^{-1})(z_1q-z_3q^{})(z_1-z_2q^{-2})} {(z_1-z_3)(z_3-z_1q^{-2})}
X^-(z_1)X^-(z_2)X^-(z_3)- $$
$$ \frac{((z_1-z_3q^{-2})(z_1-z_3q^{-1})(z_1q-z_3q^{})(z_2-z_1q^{-2})(z_2-z_1q^{})}
 {(z_1-z_3)(z_3-z_1q^{-2})(z_1-z_2q^{})}
 X^-(z_2)X^-(z_1)X^-(z_3) )= 0, $$
where the coefficient functions of the
relations above are expanded in the region 
 the expansion region 
of the corresponding monomial of the product of $X^{\pm}(z_j)$.
\end{prop} 
 
The proof is a simple calculation. (See also \cite{Er})
It is not very difficult 
to show that this relation will give us the classical Serre relations, but 
unfortunately, we still do not know how to write down a simple Serre relation 
like that of $U_q(\hat {\frak sl}(3))$.

\begin{defn}
$U_q(g, \bar {\frak s})$ is the quotient algebra of $U_q(g, {\frak s})$ with 
the kernel defined by the q-Serre relations,  
$\bar a(m)=0, m>0$, $\bar b(m), m<0$ and $a(0)=-b(0)$, and  $g(z)$ is expanded around 0.
\end{defn}

\begin{thm} 
$U_q(g, \bar {\frak s})$ is also a Hopf algebra, whose Hopf algebra structure is the 
same as that of $U_q(g, {\frak s})$. 
\end{thm} 
This can be proven by calculation as in \cite{DI}

Another immediately result we can derive is that $U_q(g, {\frak s})$ \cite{Di} \cite{GZ}
can be identified with $U_q(\hat {\frak osp}(1,2))$ , 
where $U_q(\hat {\frak osp}(1,2))$ is derived from the FRTS realization using R-matrix 
and L-operators.

\begin{thm}
$U_q(g ,{\frak s})$ with the q-Serre relation is isomorphic to 
$U_q(\hat {\frak osp}(1,2)$. 
\end{thm}

It is still an open and interesting question to study 
$U_q(g, {\frak s})$ given by other function $g(z)$.

\begin{ack}
We  would like to thank S. Khoroshkin and   T. Hodges for their  advice and encouragement. 
 Boris Feigin would like to thank the support form the grants: RFBR 99-01-01169 and  INTAS-OPEN-97-1312. \end{ack}

\end{document}